# Elementary evaluations of some Euler sums

Donal F. Connon

**Abstract**

This short note contains elementary evaluations of some Euler sums.

It was shown by the author in [1, Eq. (4.4.155zi)] that

$$(-1)^{p+1} n \int_0^1 (1-t)^{n-1} \log^p t \, dt = p! \sum_{k=1}^n \binom{n}{k} \frac{(-1)^k}{k^p}$$

and with $p = 1$ we obtain

$$n \int_0^1 (1-t)^{n-1} \log t \, dt = \sum_{k=1}^n \binom{n}{k} \frac{(-1)^k}{k} = -H_n^{(1)}$$

where $H_n^{(r)}$ are thee generalised harmonic numbers $H_n^{(r)} = \sum_{k=1}^n \frac{1}{k^r}$.

The summation then gives us

$$\sum_{n=1}^\infty \frac{H_n^{(1)}}{n^q} = -\sum_{n=1}^\infty \frac{1}{n^{q-1}} \int_0^1 (1-t)^{n-1} \log t \, dt$$

$$= -\int_0^1 \sum_{n=1}^\infty \frac{(1-t)^{n-1}}{n^{q-1}} \log t \, dt$$

and we obtain

$$\sum_{n=1}^\infty \frac{H_n^{(1)}}{n^q} = -\int_0^1 \frac{Li_{q-1}(1-t) \log t}{1-t} \, dt$$

where $Li_s(x)$ are the polylogarithm functions $Li_s(x) = \sum_{n=1}^\infty \frac{x^n}{n^s}$.

With $q = 2$ we have

$$\sum_{n=1}^{\infty}\frac{H_n^{(1)}}{n^2} = -\int_0^1 \frac{Li_1(1-t)\log t}{1-t}dt = \int_0^1 \frac{\log^2 t}{1-t}dt$$

With integration by parts we have

$$\int \frac{\log^2 t}{1-t}dt = 2Li_3(t) - 2Li_2(t)\log t - \log(1-t)\log^2 t + c$$

and we end up with the well-known result originally derived by Euler

$$\sum_{n=1}^{\infty}\frac{H_n^{(1)}}{n^2} = 2\varsigma(3)$$

With $q=3$ we have

$$\sum_{n=1}^{\infty}\frac{H_n^{(1)}}{n^3} = -\int_0^1 \frac{Li_2(1-t)\log t}{1-t}dt$$

and we see that

$$\int \frac{Li_2(1-t)\log t}{1-t}dt = \frac{1}{2}\left[Li_2(1-t)\right]^2 + c$$

This then gives us the known result

$$\sum_{n=1}^{\infty}\frac{H_n^{(1)}}{n^3} = \frac{1}{2}\varsigma^2(2)$$

Georghiou and Philippou [3] derived the following formula in 1983

$$\sum_{n=1}^{\infty}\frac{H_n^{(1)}}{n^{2p+1}} = \frac{1}{2}\sum_{j=2}^{2p}(-1)^j \varsigma(j)\varsigma(2p-j+2)$$

and we may therefore deduce the corresponding integral

$$\int_0^1 \frac{Li_{2p}(1-t)\log t}{1-t}dt = \frac{1}{2}\sum_{j=2}^{2p}(-1)^j \varsigma(j)\varsigma(2p-j+2)$$

We also have the product

$$\left[H_n^{(1)}\right]^2 = n^2 \int_0^1 (1-t)^{n-1} \log t \, dt \int_0^1 (1-u)^{n-1} \log u \, du$$

and we get the summation

$$\sum_{n=1}^\infty \frac{\left[H_n^{(1)}\right]^2}{n^q} = \sum_{n=1}^\infty \frac{1}{n^{q-2}} \int_0^1 (1-t)^{n-1} \log t \, dt \int_0^1 (1-u)^{n-1} \log u \, du$$

$$= \int_0^1 \int_0^1 \frac{Li_{q-2}[(1-t)(1-u)] \log t \log u \, dudt}{(1-t)(1-u)}$$

Letting $q = 2$ we have

$$\sum_{n=1}^\infty \frac{\left[H_n^{(1)}\right]^2}{n^2} = \int_0^1 \int_0^1 \frac{Li_0[(1-t)(1-u)] \log t \log u \, dudt}{(1-t)(1-u)}$$

$$= \int_0^1 \int_0^1 \frac{\log t \log u \, dudt}{1-(1-t)(1-u)}$$

The Wolfram Integrator instantly gives us

$$\int \frac{\log t}{1-(1-t)(1-u)} dt = \frac{1}{1-u} \log t \log\left(\frac{1-(1-t)(1-u)}{u}\right) + \frac{1}{1-u} Li_2\left(-\frac{(1-u)t}{u}\right) + c$$

and hence we have

$$\int_0^1 \frac{\log t}{1-(1-t)(1-u)} dt = \frac{1}{1-u} Li_2\left(-\frac{(1-u)}{u}\right)$$

We then obtain

$$\sum_{n=1}^\infty \frac{\left[H_n^{(1)}\right]^2}{n^2} = \int_0^1 \frac{\log u}{1-u} Li_2\left(-\frac{(1-u)}{u}\right) du$$

and it may be noted that the Wolfram Integrator cannot evaluate this integral.

We now employ the identity originally obtained by Landen [4] in 1780

$$Li_2\left(-\frac{(1-u)}{u}\right) = -\frac{1}{2}\log^2 u - Li_2(1-u)$$

with the result that

$$\int_0^1 \frac{\log u}{1-u} Li_2\left(-\frac{(1-u)}{u}\right) du = -\frac{1}{2}\int_0^1 \frac{\log^3 u}{1-u} du - \int_0^1 \frac{Li_2(1-u)\log u}{1-u} du$$

It is well known that $\int_0^1 \frac{\log^3 u}{1-u} du = -6\varsigma(4)$ and the second integral has already been evaluated above. We accordingly get

$$\int_0^1 \frac{\log u}{1-u} Li_2\left(-\frac{(1-u)}{u}\right) du = \frac{17}{4}\varsigma(4)$$

and hence we have

$$\sum_{n=1}^\infty \frac{\left[H_n^{(1)}\right]^2}{n^2} = \frac{17}{4}\varsigma(4)$$

as originally discovered by de Doelder [2] in 1991.

Letting $q = 3$ we have

$$\sum_{n=1}^\infty \frac{\left[H_n^{(1)}\right]^2}{n^3} = \int_0^1\int_0^1 \frac{Li_1[(1-t)(1-u)]\log t \log u\, dudt}{(1-t)(1-u)}$$

$$= -\int_0^1\int_0^1 \frac{\log[1-(1-t)(1-u)]\log t \log u\, dudt}{(1-t)(1-u)}$$

The Wolfram Integrator can evaluate the integral

$$\int \frac{\log[1-(1-t)(1-u)]\log t}{(1-t)(1-u)} dt$$

but additional work is required to check if the definite integral exists.